\documentclass[10pt]{article}
\newcommand{\documentdate}{12 I 2026}
\usepackage{a4wide,latexsym,graphicx,amsmath,amssymb,varioref}
\usepackage{graphics,xcolor,dsfont}

\newcommand{\al}[1]{{\footnotesize{\sf #1}}}
\title{Iteration complexity of the Difference-of-Convex Algorithm for
  unconstrained optimization: a simple proof}

\author{
S. Gratton\thanks{Universit\'{e} de Toulouse, INP, IRIT, Toulouse, France. Email:
     serge.gratton@enseeiht.fr. Work partially supported by 3IA Artificial and
     Natural Intelligence Toulouse Institute (ANITI), French "Investing for the Future
     - PIA3" program under the Grant agreement ANR-19-PI3A-0004"}
~and Ph. L. Toint\thanks{NAXYS, University of Namur, Namur, Belgium. Email: philippe.toint@unamur.be}
}

\DeclareMathAlphabet{\pazocal}{OMS}{zplm}{m}{n}

\newcommand{\calO}{{\pazocal{O}}}

\newcommand{\beqn}[1]{\begin{equation}\label{#1}}
\newcommand{\eeqn}{\end{equation}}
\newcommand{\req}[1]{(\ref{#1})}
\newcommand{\ms}{\;\;\;\;}
\newcommand{\tim}[1]{\;\; \mbox{#1} \;\;}

\newtheorem{theorem}{Theorem}[section]
\newtheorem{lemma}[theorem]{Lemma}

\newcommand{\numsection}[1]{\section{#1}\setcounter{equation}{0}}
\newtheorem{corollary}{Corollary}

\renewcommand{\thefootnote}{(\arabic{footnote})}
\newcounter{algo}[section]
\renewcommand{\thealgo}{\thesection.\arabic{algo}}

\newcommand{\algo}[3]{\refstepcounter{algo}
\begin{center}\begin{figure}[htbp]
\framebox[\textwidth]{
\parbox{0.95\textwidth} {\vspace{\topsep}
{\bf Algorithm \thealgo : #2}\label{#1}\\
\vspace*{-\topsep} \mbox{ }\\
{#3} \vspace{\topsep} }}
\end{figure}\end{center}}
\newcommand{\bpr}{{\bf Proof.} \hspace{1.5mm}}
\newcommand{\epr}{\hfill $\Box$ \vspace*{1em}}
\newcommand{\proof}[1]{
\begin{list}{}{
\setlength{\topsep}{0.0pt}
\setlength{\partopsep}{0.0pt}
\setlength{\leftmargin}{0.025\textwidth}
\setlength{\rightmargin}{0.5\leftmargin}
\setlength{\labelwidth}{0.5\leftmargin}
\setlength{\labelsep}{0.25\leftmargin}}
\item \bpr #1 \epr \noindent
\end{list}}
\newcommand{\lthm}[2]{\vspace{\baselineskip} 
\noindent\framebox[\textwidth]{\parbox{0.95\textwidth}{
\begin{theorem} \label{#1} \rm #2 \end{theorem} } } \vspace{\baselineskip} }

\renewcommand{\Re}{\hbox{I\hskip -2pt R}}
\newcommand{\smallRe}{\hbox{\footnotesize I\hskip -2pt R}}

\newcommand{\sfrac}[2]{{\scriptstyle \frac{#1}{#2}}}

\newcommand{\eqdef}{\stackrel{\rm def}{=}}
\newcommand{\tal}[1]{{\normalsize {\sf #1}}}
\newcommand{\half}{\sfrac{1}{2}}
\newcommand{\flow}{f_{\rm low}}

\newcommand{\comment}[1]{}
\DeclareMathOperator*{\argmin}{argmin}

\newcommand{\ckh}{\lceil k/2 \rceil}
\newcommand{\fkh}{\lfloor k/2 \rfloor}

\date{\documentdate}
\begin{document}

\renewcommand{\thefootnote}{\fnsymbol{footnote}}
\maketitle
\renewcommand{\thefootnote}{\arabic{footnote}}

\begin{abstract}
We propose a simple proof of the worst-case iteration complexity for
the  Difference of Convex functions Algorithm \al{DCA} for
unconstrained minimization, showing that the global rate of
convergence of the norm of the objective function's gradients at the
iterates converge to zero like $o(1/k)$. A small example is also
provided indicating that the rate cannot be improved.
\end{abstract}

{\small
 \al{DCA} algorithm, iteration complexity.
}

\numsection{Introduction}

Proposed by Pham Dinh and Souad in \cite{PhamSoua86} following a
first decomposition approach by Mine and Fukushima \cite{MineFuku81},
the possibility to decompose a possibly nonconvex function into the
difference of two convex functions (DC) has raised a considerable
interest, both practically and theoretically, as is testified by the
many references on this topic (see, for instance,
\cite{AragFlemVuon18,GaudGialMiglBagi18,LeTh00,LeThDinh05,PhamLeTh97,
  PhamHoai98,LeThPhanDinh21,Niu22}
or the survey \cite{LeThPham23}).  This interest is often motivated by
the fact that functions of this particular class exhibit some
structure, which one hopes to exploit in order to improve computational
efficiency.  It is therefore natural to wonder if
this hope of improved performance can be translated into an improved
worst-case complexity.  As it turns out, this question was solved in
\cite{AbbadeKlZama24} and further refined in \cite{RotaPanaGlin25},
using the complex machinery of performance analysis as proposed in
\cite{DrorTebo14}.  The proofs in these papers are non-trivial as this
framework entails the semi-definite programming relaxation of an infinite
dimensional optimization problem with an infinite number of
constraints.  The purpose of this short note is to show that both
lower and upper tight complexity bounds may be obtained by a much
simpler, direct and explicit approach, which has the advantage of
conciseness but also of improved interpretability.

Thus, in what follows, we (re-)analyze the global convergence rate of the
iterations of the standard algorithm for unconstrained minimization of
DC functions and show that it is in order identical to that of the
steepest descent method applied on the unstructured problem.

More specifically, we consider solving the smooth unconstrained
problem 
\beqn{problem}
\min_{x\in\smallRe^n} f(x) = g(x) - h(x)
\eeqn 
where $g$ and $h$ are continuously differentiable convex functions from
$\Re^n$ into $\Re$, the gradient of $h$ being Lipschitz continuous.
To do so, we will use the standard \al{DCA} (DC algorithm) \cite{MineFuku81,PhamSoua86}, which
is stated \vpageref{DCA}.

\algo{DCA}{\tal{DCA}}{
  \begin{description}
  \item[Step 0: Initialization.]  A starting point $x_0$ is
    given. Set $k=0$.
  \item[Step 1: Step computation.]   For $k \geq 0$, set
    \beqn{xkp1}
    x_{k+1} = \argmin_{x\in \Re^n} g(x) - \nabla_x h(x_k)^T(x-x_k)
   \eeqn
 \end{description}  
}

Observe that the function $g(x) - \nabla_x h(x_k)^T(x-x_k)$ (used to
define the next iterate in Step~1) inherits
the convexity of $g$. As a consequence the solution of the global
minimization in \req{xkp1} can make use of the arsenal of methods for
convex problems, whose efficiency often vastly exceeds that of methods
for general nonconvex ones. While this does influence the total
evaluation complexity of the \al{DCA}, as we briefly discuss
below, we focus in this note on its iteration complexity, or,
equivalently, on the global rate of convergence of its iterations.

We also immediately note that, because $h$ is convex,
$g(x) - h(x_k) - \nabla_x h(x_k)^T(x-x_k)$ is an overestimate of
$f(x)$ and therefore 
\beqn{decr}
f(x_{k+1}) \le f(x_k)
\eeqn
for all $k\ge0$.

\numsection{Global rate of convergence}

We start our analysis of the global rate of convergence of the
algorithm by establishing an upper bound on the global rate of
decrease of the sequence $\{\|\nabla_x f(x_k)\|\}$ to zero
(see \cite[Theorem~3.1]{AbbadeKlZama24}).

\lthm{upper-complexity}{
  Consider applying the \al{DCA} to
  problem \req{problem}. Suppose that $g$ is continuously
  differentiable and strongly convex with constant $\mu>0$, that is 
  \beqn{strongly-convex}
  g(y) \geq g(x) + \nabla_xg(x)^T(y-x) + \mu\|y-x\|^2 \tim{ for all } x,y \in \Re^n,
  \eeqn
  and that $h$ is convex, continuously differentiable and that its gradient is
  Lipschitz continuous with constant $L_h$, that is
  \beqn{Lipschitz}
  \|\nabla_xh(x)-\nabla_xh(y)\| \le L_h \|x-y\| \tim{ for all } x,y \in \Re^n.
  \eeqn
  Suppose finally that there exists $\flow\in \Re$ such that
  $f(x)\ge\flow$ for all $x\in \Re^n$.
  Then, for $k\ge 2$,
  \beqn{the-rate}
  \frac{1}{\fkh} \sum_{i=\ckh}^k \|\nabla_x f(x_i)\|^2
  \le \frac{2L_h^2(f(x_{\ckh})-f(x_{k+1}))}{\mu\, k}
  = o\left(\frac{1}{k}\right)
  \eeqn
}

\proof{
  The strong convexity of $g$ implies that
  \beqn{e1}
  g(x_k)\ge g(x_{k+1}) + \nabla_x g(x_{k+1})^T(x_k-x_{k+1}) +\mu \|x_{k+1}-x_k\|^2.
  \eeqn
  while the convexity of $h$ gives that
  \beqn{e2}
  h(x_{k+1})\ge h(x_k) + \nabla_x h(x_k)^T(x_{k+1}-x_k).
  \eeqn
  Observe now that \req{xkp1} implies that
  \beqn{e3}
  \nabla_xg(x_{k+1}) = \nabla_x h(x_k).
  \eeqn
  Thus \req{e1} gives that
  \beqn{e4}
  g(x_k)\ge g(x_{k+1}) + \nabla_x h(x_k)^T(x_k-x_{k+1}) +\mu
  \|x_{k+1}-x_k\|^2.
  \eeqn
  Summing \req{e2} and \req{e4} then yields that
  \[
  g(x_k)-h(x_k) \ge g(x_{k+1})-h(x_{k+1}) + \mu\|x_{k+1}-x_k\|^2.
  \]
  Summing now this inequality over the last $\ckh$ iterations gives that
  \beqn{e5}
  f(x_{\ckh})-f(x_{k+1})
  = \sum_{i=\ckh}^k  \left[\big(g(x_{i+1})-h(x_{i+1})\big)-\big(g(x_{i})-h(x_{i})\big)\right]
  \ge  \mu \sum_{i=\ckh}^k \|x_{i+1}-x_i\|^2.
  \eeqn
  Now, using \req{e3} and \req{Lipschitz}, we see that, for all $i\geq0$,
  \[
  \|\nabla_x f(x_i)\|
  = \|\nabla_x g(x_i)-\nabla_x h(x_i)\|
  =  \|\nabla_x h(x_{i+1})-\nabla_x h(x_i)\|
  \le L_h\|x_{i+1}-x_i\|.
  \]
  This and \req{e5} then imply that
  \[
  \sum_{i=\ckh}^k\|\nabla_x f(x_i)\|^2
  \le L_h^2\sum_{i=\ckh}^k\|x_{i+1}-x_i\|^2
  \le \frac{L_h^2(f(x_{\ckh}) - f(x_{k+1}))}{\mu}.
  \]
  Thus, dividing by $\fkh+1$, 
  \beqn{ef}
  \frac{1}{\fkh+1} \sum_{i=\ckh}^k \|\nabla_x f(x_i)\|^2
  \le \frac{L_h^2(f(x_{\ckh})-f(x_{k+1}))}{\mu(\fkh+1)}
  \le \frac{2L_h^2(f(x_{\ckh})-f(x_{k+1}))}{\mu\, k}.
  \eeqn
  Since the sequence $\{f(x_k)\}$ is monotonically decreasing because
  of \req{decr} and bounded below by assumption, it must be convergent
  and thus
  \[
  \lim_{k\rightarrow \infty} (f(x_{\ckh}) - f(x_{k+1})) = 0.
  \]
  Thus \req{ef} yields \req{the-rate}.
} %epr

\noindent
Note that inequality \req{e5} recovers
a variant of a known result on the sequence of iterates
(\cite[Proposition~2]{PhamLeTh97}, \cite[Corollary~1]{LeThPhanDinh21},
\cite[Theorem~1.2]{AbbadeKlZama24}).

Should one wish to stop the algorithm as soon as an iterate $x_k$ is
found such that
\beqn{eps-termination}
\|\nabla_xf(x_k)\| \le \epsilon
\eeqn
for some user-prescribed tolerance $\epsilon >0$,
Theorem~\ref{upper-complexity} ensures that at most $o(\epsilon^{-2})$
iterations are needed for termination.  While this bound is
theoretically better that $\calO(\epsilon^{-2})$ of
\cite[Theorem~3.1]{AbbadeKlZama24}, we now show that the number of
iterations necessary to achieve \req{eps-termination} may be
arbitrarily close to $\calO(\epsilon^{-2})$, providing a lower 
bound on iteration complexity for the \al{DCA}. This is achieved by
providing an explicit example of ``slow'' convergence, much as in
\cite[Section~3]{AbbadeKlZama24}, but using simpler and more
interpretable functions.

\lthm{lower-complexity}{
Let the univariate function g be defined by $g(x) = \half x^2$. Then,
for any $\delta >0$,
there exists a convex continuously differentiable function $h$ from
$\Re$ to $\Re$ with Lipschitz continuous gradient such that the
sequence of iterates of the \al{DCA} applied to
problem \req{problem} starting from $x_0=0$ is such that 
\beqn{slow}
\|\nabla_x f(x_k) \|^2 =  \frac{1}{(k+1)^{1+2\delta}}.
\eeqn
}

\proof{
  Define the sequence $\{x_k\}$ by
  \beqn{xk}
  x_{k+1} = x_k - \frac{1}{(k+1)^{\half+\delta}} \ms (k\ge 0)
  \eeqn
  which implies that $\{x_k\}$ tends to $-\infty$. We now construct a
  function $h$ by first imposing the conditions
  \beqn{sync}
  \nabla_xh(x_k) = x_{k+1}\ms (k\ge 0)
  \eeqn
  and then defining $h$ on $[x_{k+1},0]$ as the continuously differentiable
  piecewise quadratic interpolating the conditions \req{sync} and such
  that $h(x_0) = h(0)= 0$. One easily checks that its second derivative on
  the interval $[x_{k+1},x_k]$ is given by
  \[
  \frac{\nabla_xh(x_{k+1})-\nabla_xh(x_k)}{x_{k+1}-x_k}
  = \frac{x_{k+2}-x_{k+1}}{x_{k+1}-x_k}
  = \frac{(k+1)^{\half+\delta}}{(k+2)^{\half+\delta}}
  \in (0,1]
  \]
  and we then define
  \beqn{hk-def}
  \begin{aligned}
  h(x_{k+1})
  &= h(x_k) +\nabla_xh(x_k)(x_{k+1}-x_k)+\half \nabla_x^2h(x_k)(x_{k+1}-x_k)^2\\
  &= h(x_k) - \frac{x_{k+1}}{(k+1)^{\half+\delta}}
            + \frac{1}{2(k+1)^{1+2\delta}}\left(\frac{k+1}{k+2}\right)^{\half+\delta}
  \end{aligned}
  \eeqn
  Hence $h$ is convex on $(-\infty,0]$ and, using
  \cite[Theorem~A.8.5]{CartGoulToin22}, its gradient is Lipschitz
  continuous with constant equal to one. It is then easy to prolongate
  $h(x)$ for $x>0$ by setting
  \[
  h(x) = h(0)+\nabla_xh(0)\, x = -x
  \]
  for $x>0$ without altering its convexity or the Lipschitz
  continuity of its gradient. Writing $g(x)$ as its second-order
  Taylor's series, we also have that, for $k\ge 0$,
  \[
  g(x_{k+1}) = g(x_k) - \frac{x_k}{(k+1)^{\half+\delta}}+\frac{1}{2(k+1)^{1+2\delta}},
  \]
  and thus, combining this relation with \req{hk-def}, that
  \[
  f(x_{k+1})
  = f(x_k) - \frac{x_k-x_{k+1}}{(k+1)^{\half+\delta}}
    +\frac{1}{2(k+1)^{1+2\delta}}
    \left[1-\left(\frac{k+1}{k+2}\right)^{\half+\delta}\right]
  \ge f(x_k) - \frac{1}{(k+1)^{1+2\delta}}.
  \]
  Therefore, for all $k\geq0$,
  \[
  f(x_{k+1}) \ge f(x_0) - \sum_{i=0}^k \frac{1}{(i+1)^{1+2\delta}}
  \ge - \zeta(1+2\delta),
  \]
  where $\zeta(\cdot)$ is the Riemann zeta function. Moreover, since
  both $g(x)$ and $h(x)$ have Lipschitz continuous gradients with
  unit constant, we have that, for $x\in [x_{k+1},x_k]$,
  \[
  f(x) \ge f(x_k) - \frac{1}{(k+1)^{1+2\delta}} - \|x-x_k\|^2
  \ge f(x_k) - \frac{2}{(k+1)^{1+2\delta}}
  \ge f(x_k) - 2,
  \]
  where we have used \req{xk} to deduce the last inequality.
  This in turn implies that, for all $x \in (-\infty, 0]$,
  \beqn{fbbelow}
  f(x)\ge -\zeta(1+2\delta) -2 =\eqdef \flow.
  \eeqn
  Since $f(x) = \half x^2 + x >0 $ for $x>0$, we see that $g$ and $h$ satisfy all the conditions
  required by Theorem~\ref{upper-complexity}.  We may then apply the
  \al{DCA}  to minimize $f=g-h$. The definition of its iteration
  in \req{xkp1} gives that 
  \beqn{iterok}
  \argmin_{x\in\smallRe} \half x^2 - x_{k+1}x = x_{k+1}.
  \eeqn
  and the iterates generated by the algorithm therefore coincide with
  the sequence $\{x_k\}$ defined in \req{xk}.  Moreover, we deduce
  from \req{xk} and \req{sync} that
  \[
  \|\nabla_x f(x_k)\| = \|x_k - x_{k+1}\| = \frac{1}{(k+1)^{\half+\delta}},
  \]
  which is \req{slow}.
} % epr

\begin{figure}
  \begin{center}
    \includegraphics[width=8cm]{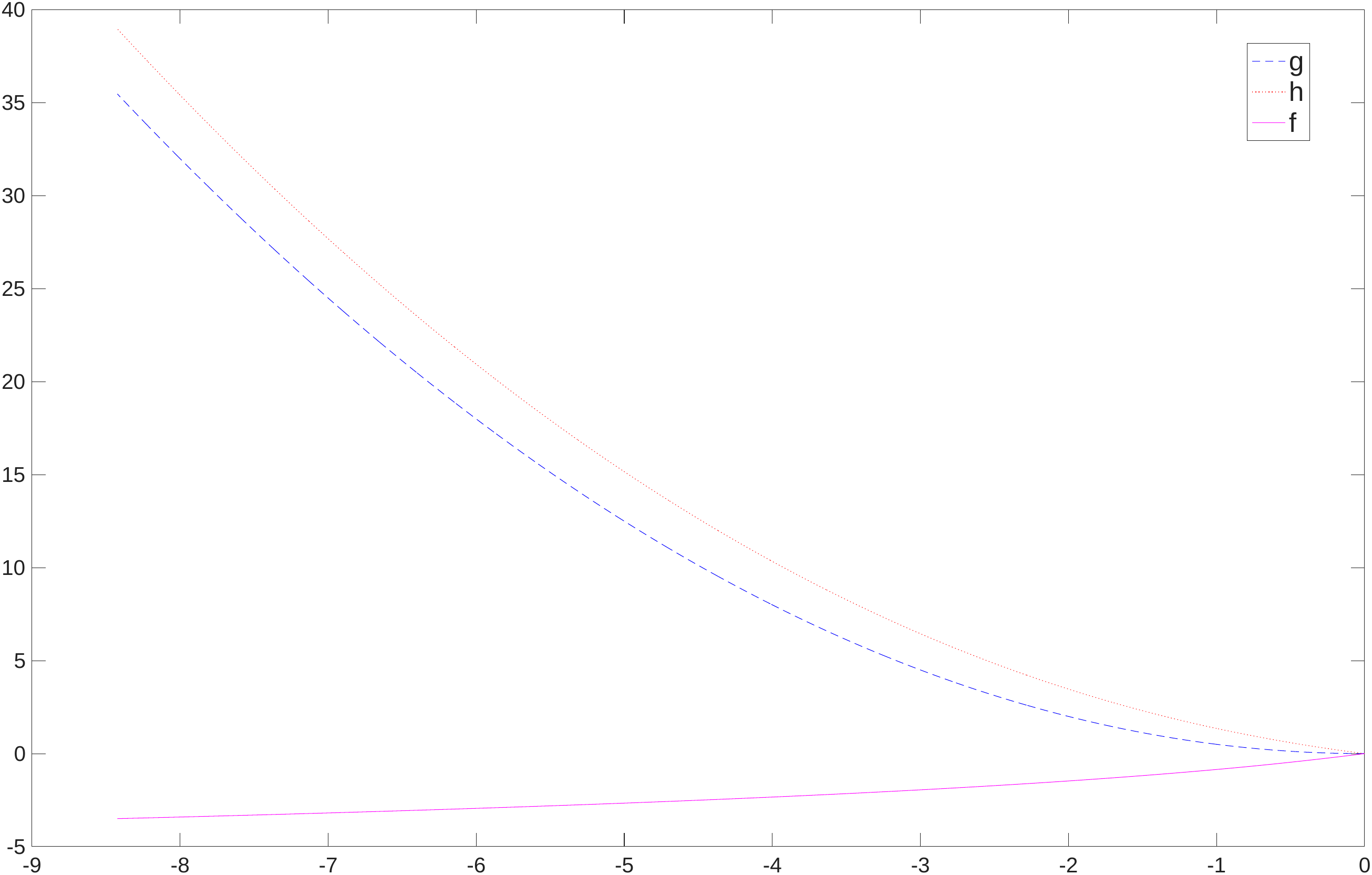}
  \caption{\label{fig:fgh}The shapes of $f$, $g$ and $h$ on
    $[x_{25},x_0]$}
  \end{center}
\end{figure}

\noindent
The shapes of $f$, $g$ and $h$ on $[x_{25},x_0]$ are shown in
Figure~\ref{fig:fgh}.

\numsection{Discussion}

We have provided proofs for the lower an upper bounds on iteration
complexity for the \al{DCA}, simplifying results of
\cite{AbbadeKlZama24}. The main conclusion that can be drawn from
these results is that the minimization of DC functions using the
\al{DCA} does not result in an better iteration complexity than that
obtained for the steepest-descent method applied on $f$ and ignoring
the DC structure \cite{GratSimToin25}. Thus distinguishing DC
functions from general convex ones does not bring any benefit from the
worst-case complexity point of view, should this algorithm be used. It
is however important to remember that this does not mean that the
practical perfoemance of the \al{DCA} is as slow as that of steepest
descent (after all, this is also the case of several efficient
variants of Newton's method, as shown in
\cite[Theorem~3.1.1]{CartGoulToin22}).

Our example of Theorem~\ref{lower-complexity} also shows that using
the \al{DCA} to exploit structure might be detrimental from the point
of view of complexity. Indeed, it is easy to verify that the Hessian
of $f$ is strictly positive on each interval $[x_{k+1},x_k]$, and the
function $f$ is therefore convex (albeit not strongly convex). The
performance of the \al{DCA} on that convex function is therefore worse
than what would be achieved by a good algorithm for general convex
function (see \cite[Chapter~2]{Nest18} for an extensive discussion).

We also note that evaluation complexity (that is the maximum number of
oracle\footnote{objective function's value and/or gradient.}
evaluations needed to achieve \req{eps-termination}) is bounded below by
iteration complexity, because the solution of the inner minimization
subproblem of \req{xkp1} does require at least one oracle evaluation.

We finally observe that our unidimensional example of slow convergence
presented in Theorem~\ref{lower-complexity} can also
straightforwardly be extended to any dimension by considering
functions of an arbitrary number of variables but whose value and
gradient only depend on one. It is also remarkable that the example
is independent of any
termination rule, at variance with examples of slow convergence in
\cite{CartGoulToin22}.  This is due to the exploitation of the proof
technique of \cite{GratSimToin25} in
Theorem~\ref{lower-complexity}. Indeed, our example does not extend to
the case where $\delta=0$ unless a termination rule of the type
\req{eps-termination} is introduced, because the constructed
function $f$ is no longer bounded below in that case.

{\footnotesize

\section*{\footnotesize Acknowledgement}

Philippe Toint is grateful for the continued and friendly support of
the APO team at Toulouse IRIT. Both authors thank M. Al-Baali (Sultan
Qaboos University, Oman) for organizing the NAOVI-2026 conference in
Muscat, during which the question of the complexity of DC optimization
emerged after a talk by A. Bagirov.
  
\bibliography{/home/philippe/bibs/refs}
\bibliographystyle{plain}

}% footnotesize

\end{document}